\documentclass{amsart} 

\newtheorem{theorem}{Theorem}[section]
\newtheorem{lemma}[theorem]{Lemma}
\newtheorem{corollary}[theorem]{Corollary}
\newtheorem{proposition}[theorem]{Proposition}

\theoremstyle{definition}

\newtheorem{example}[theorem]{Example}
\newtheorem{remark/example}[theorem]{Remark/Example}

\numberwithin{equation}{section}

\def\PP{ {\bf P} }

\def\CC{ {\bf C} }

\newcommand{\GL}{\operatorname{GL}}
\newcommand{\rank}{\operatorname{rank}}
\newcommand{\codim}{\operatorname{codim}}
\numberwithin{equation}{section}
 
\begin{document}

\title{Gr\"obner bases for spaces of quadrics of codimension $3$ }
\author{Aldo Conca }
\address{ Dipartimento di Matematica, Universit\'a di Genova, Via Dodecaneso 35, I-16146 Genova, Italia }
\email{conca@dima.unige.it }
%\date{}

\begin{abstract}
Let $R=\oplus_{i\geq 0} R_i$ be an Artinian standard graded $K$-algebra defined by quadrics. 
Assume that $\dim R_2\leq 3$ and that $K$ is algebraically closed of characteristic $\neq 2$. We show that $R$ is defined by a Gr\"obner basis of quadrics with, essentially, one exception. The exception is given by $K[x,y,z]/I$ where $I$ is a complete intersection of $3$ quadrics not containing a square of a linear form. 
\end{abstract} 
\maketitle

\section{introduction} 
A standard graded $K$-algebra $R$ is an algebra of the form
$R=K[x_1,\dots,x_n]/I$ where $K[x_1,\dots,x_n]$ is a polynomial ring over
the field $K$ and $I$ is a homogeneous ideal with respect to the grading
$\deg(x_i)=1$. The algebra $R$ is said to be { quadratic } if $I$ is
generated by quadrics (i.e.~homogeneous elements of degree two) and $R$ is
said to be { Koszul } if $K$ admits a free resolution as an $R$-module 
whose maps are given by matrices of linear forms. We say that an algebra $R$ is { 
G-quadratic} if its defining ideal has a Gr\"obner basis of quadrics with respect to some system of coordinates and some term order. G-quadratic algebras
are Koszul and Koszul algebras are quadratic. Neither of these implications can be reversed in general, see \cite{F}. 

Given a graded $K$-algebra $R$ we can consider its trivial fiber extension $R'=R\circ K[y]/(y)^2$ where $y=y_1,\dots,y_m$ is a set of variables. Here $\circ$ denotes the fibre product of $K$-algebras.  It is known that the properties  of  being quadratic, Koszul, G-quadratic, as well as $\dim R_i$ for $i>1$ are unaffected by trivial fiber extensions, see \cite[Lemma 4]{C}, \cite[Theorem 4]{BF}.

Backelin showed in \cite[4.8]{B} that if $R=\oplus_{i\geq
0}R_i$ is a quadratic standard graded $K$-algebra with $\dim R_2\leq 2$ 
then $R$ is Koszul. We have shown in \cite{C} that, under the same assumptions, $R$ is $G$-quadratic with  only one exception (up to trivial fiber extensions and changes of coordinates) given by the $K$-algebra $K[x,y,z]/(x^2,xy,y^2-xz,yz)$. 

The goal of this paper is to prove the following theorem. 

\begin{theorem} 
\label{main}
Let $K$ be an algebraically closed field of characteristic $\neq 2$. Let $R$ be a standard graded $K$-algebra which is quadratic, Artinian and with $\dim R_2=3$. Then: 
\begin{itemize}
\item[(1)] $R$ is Koszul, $R_i=0$ for $i>3$ and $\dim R_3\leq 1$. Furthermore $\dim R_3=1$ if and only if $R$ is a trivial fiber extension of $K[x,y,z]/I$ where $I$ is a complete intersection of $3$ quadrics. 
\item[(2)] $R$ is G-quadratic iff  it is not a trivial fiber extension of $K[x,y,z]/I$ where $I$ is a complete intersection of $3$ quadrics not containing a square of a linear form.  In particular, if $R_3=0$, then $R$ is G-quadratic. 
\end{itemize} 
\end{theorem} 

We obtain the following corollaries. 
 
 \begin{corollary} 
 \label{CM-cor}
Let $R$ be a quadratic Cohen-Macaulay standard graded $K$-algebra. Denote by $(h_0,h_1,h_2,\dots)$ its  $h$-vector  and assume that $h_2=3$. Then: 
\begin{itemize} 
\item[(1)] $R$ is Koszul, $h_i=0$ for every $i>3$ and $h_3\leq 1$. Furthermore $h_3=1$ if and only if  the degree $1$ component of the socle of $R$ has dimension $h_1-3$.  
\item[(2)] If $h_3=0$  then $R$ is G-quadratic.
\end{itemize} 
\end{corollary}

 \begin{corollary} 
 \label{CM-cor2}
Let $R$ be a quadratic Cohen-Macaulay algebra. If   $e(R)\leq \codim(R)+4$ then $R$ is Koszul. 
\end{corollary}

In \ref{CM-cor2} $e(R)$ denotes the degree or multiplicity of $R$  and $\codim(R)$ its codimension.  To see that  \ref{CM-cor} follows form \ref{main} one just considers an  Artinian reduction of $R$.  But then  \ref{CM-cor2} follows from \ref{CM-cor} in combination with \cite[Corollary 9]{C}.   In particular a (non-degenerate) set of at most $n+4$ points of $\PP^n$ defiend by quadrics is Koszul, a special case of a recent conjecture of Polishchuk \cite{P}.  Note also that there are non-Koszul quadratic Cohen-Macaulay algebras (even domains) with $e(R)=\codim(R)+5$, see \cite[Sect.4]{CRV}.

For standard facts on Gr\"obner bases we refer the reader to \cite{KR}. 
The results and the examples presented   were discovered  by extensive computer algebra experiments performed with CoCoA \cite{Co}.

\section{Proof of the main result} 

Let $K$ be an algebraically closed field of characteristic not $2$. Let $R$ be a standard graded $K$-algebra.  The rank of $x\in R_1$, denoted by $\rank x$, is by definition $\dim xR_1$. Note that $\rank x=0$  for some non-zero $x\in R_1$ iff $R$ is a trivial fiber extension. 

As we said already, the properties under discussion, being quadratic, Koszul, G-quadratic, as well as $\dim R_i$ for $i>1$ are unaffected by trivial fiber extensions, see \cite[Lemma 4]{C}, \cite[Theorem 4]{BF}.  Hence in the proof of   \ref{main} we may assume that  $\rank x>0$ for every the non-zero $x\in R_1$. 
  
The case $n=3$  is easy: $R$ is a complete intersection of $3$ quadrics and \ref{main} is proved in \cite[Sect.6.1]{CRV}.  It remains to prove: 

\begin{proposition}
\label{TH}  With  the assumption of \ref{main}, assume further that $\rank x>0$ for every non-zero $x\in R_1$ and that $n>3$.  It follows that  $R$ is G-quadratic and $R_3=0$. 
\end{proposition}

The main technical lemma is: 

\begin{lemma}
\label{ml} Under Assumption \ref{TH}, let $y\in R_1$ with $y^2=0$ and set $V=\{ u\in R_1 : uy=0\}$. If one of the following conditions hold
then $R$ is G-quadratic and $R_3=0$. 
\begin{itemize} 
\item[(1)] $\rank y=3$.

\item[(2)] $\rank y=2$ and there exists $z\in V$ such that $z^2\in yR_1$ and
$zR_1\not\subseteq yR_1$. 

\item[(3)] $\rank y=2$ and there exists $t\in R_1\setminus V$ such that $t^2\in yR_1$ and
 $tV\not\subseteq yR_1$. 

\item[(4)] $\rank y=1$.
\end{itemize} 
\end{lemma} 

In the proofs below the symbol $L_*$ denotes an homogeneous polynomial of degree $1$  and     $*f$ means a scalar multiple of the element  $f$. 

\begin{proof} In the four cases the proof is based on the same principle: choose a $K$-basis of $R_1$ in a suitable way, consider the associated presentation of $R$ as a quotient of a polynomial ring, translate the assumptions into quadratic equations, check that the given quadrics provide already enough leading terms to generate all the monomials of degree $3$. For simplicity of notation, we do not distinguish between the elements of $R_1$ and the variables of the polynomial ring that we use to present $R$. 

Case (1) is the standard situation, see \cite{C}. Complete $y$ to a basis of $R_1$, say $y,x_2,\dots,x_n$. Since $y^2=0$ and $yR_1=R_2$, we have polynomials $y^2$ and $x_ix_j-yL_*$ in the defining ideal $I$ of $R$. Then $(y)^2+(x_2,\dots,x_n)^2$ is contained in the ideal of leading terms of $I$ with respect the revlex order associated with $x_i>y$. This is enough to conclude that $R$ is G-quadratic and $R_3=0$. 

Case (2): Complete $y,z$ to a basis of $R_1$, say $y,z, x_3,\dots,x_{n}$. We have polynomials $y^2, yz, z^2-yL_*, x_ix_j-L_*y-L_*z$ in the defining ideal $I$ of $R$. It follows that $(y,z)^2+(x_3,\dots,x_n)^2$ is in the ideal of leading terms of $I$ with respect the revlex order associated with $x_i>z>y$. This is enough to conclude that $R$ is G-quadratic and $R_3=0$.

Case (3). Consider a basis $y,z_2,\dots,z_{n-2}$ of $V$ and complete it with the given $t$ and some other element $w$ to a basis of $R_1$. Use the term order revlex $y<t<z_i<w$. Polynomials of the following form are in the defining ideal: 
 $$y^2, \ \ yz_i, \ \ z_iz_j-yL_*-tL_*, \ \ t^2-yL_*, \ \ wz_i-yL_*-tL_*, \ \ w^2-yL_*-tL_*$$ 
 Set $W=\{ u\in R_1 : ut\in yR_2\}$. Note that $W$ is a space of dimension $n-1$ and contains a linear form which involves $w$ 
 otherwise we would have $tV\subseteq yR_1$ which contradict the assumption of (3). Then a polynomial of the form $wt-yL_*$ is also in the defining ideal. The leading term ideal of the defining ideal of $R$ contains 
$(y,z_2,\dots,z_{n-2})^2+(t,w)^2$. This is enough to conclude that $R$ is G-quadratic and $R_3=0$.

 Case (4). We have that $R/(y)$ is Artinian with Hilbert series 
$1+(n-1)x+2x^2+\dots$. By \cite{C} we know that there exists $t\in R_1$ such that $t^2\in yR_1$ and $R_2=yR_1+tR_1$. Complete $y$ and $t$ to a basis of $R_1$ with elements $x_3,\dots,x_n$ and use the revlex order associated with $y<x_i<t$. In the defining ideal of $R$ we have polynomials
$$y^2, \ \ x_ix_j-yL_*-tL_*, \ \ t^2-yL_*$$ 
and so initial terms $(x_3,\dots, x_n)^2+(y^2, t^2)$. Furthermore, since $yV=0$ in $R$ and $V$ has dimension $n-1$, we have inital terms $Wy$ where $W$ is a set of variables of cardinality $n-1$ containing $y$. So either all the $x_i$ are in $W$ or $t$ is in $W$. In the first case the inital term ideal contains $(x_3,\dots,x_n,y)^2+(t^2)$, in the second $(x_3,\dots,x_n)^2+(y,t)^2$. In both cases we are done. 
\end{proof}

Another auxiliary fact: 

\begin{lemma}
\label{l3} Assume $S$ is a quadratic standard graded $K$-algebra with Hilbert series $1+3x+x^2$ and $z\in S_1$ such that $z^2\neq 0$. Then there exists $u\in S_1$ such that $u^2=0$ and $uz\neq 0$. 
\end{lemma} 

\begin{proof} We argue by contratiction. Let $v\in S_1$ not a multiple of $z$. We have equations $v^2=az^2$ and
$vz=bz^2$ with $a,b\in K$. By contradiction, there is no $\alpha\in K$ such that 
 $(v+\alpha z)^2=0$ and $z(v+\alpha z)\neq 0$, that is, no $\alpha\in K$ such that 
$a+2\alpha b +\alpha^2=0$ and $b+\alpha\neq 0$. In other word, $\alpha=-b$ is the only solution of 
$a+2\alpha b +\alpha^2=0$, that is, $a=b^2$. Complete $z$ to a basis of $S_1$ with elements $t,w$. By the argument
above we have equations: 
 we have equations: 
$$t^2=b^2z^2, \ \ tz=bz^2, \ \ w^2=c^2z^2, \ \ wz=cz^2.$$ We have also an equation $wt=dz^2$. Since
$(t+w)^2=(b^2+c^2+d)z^2$ and $(t+w)z=(b+c)z^2$ 
 the argument above applied to $t+w$ yields: 
$$(b^2+c^2+2d)=(b+c)^2$$ that is $d=bc$. So the polynomials defining $S$ are 
$$t^2-b^2z^2, \ \ tz-bz^2, \ \ w^2-c^2z^2, \ \ wz-cz^2, \ \ wt-bcz^2.$$ These polynomials are
contained in the ideal $(t-bz, w-cz)$, contradicting the fact that $S$ is Artinian. 
\end{proof} 

Now we are ready to prove \ref{TH}:

\begin{proof} Fix $K$-bases of $R_1$ and $R_2$. The condition $y^2=0$ for an element $y\in R_1$ is expressed by $\dim R_2$ quadratic equations in the $\dim R_1$ coefficients of $y$. Since $\dim R_1=n>\dim R_2=3$ and $K$ is algebraically closed, there exists $y\in R_1$ non-zero such that $y^2=0$. Further $\rank y>0$ by assumption.  If $\rank y=3$
or $1$ we conclude by \ref{ml} Case (1) and Case (4). So we may assume that $\rank y=2$. Let $V=\{ u \in R_1 : uy=0\}$, $V$ is a $n-2$-dimensional subspace of $R_1$. We discuss three cases: 

\medskip 

\noindent Case 1: $V^2\not\subseteq yR_1$

\medskip 

\noindent Case 2: $V^2\subseteq yR_1$ and $VR_1\not\subseteq yR_1$ 

\medskip 

\noindent Case 3: $VR_1\subseteq yR_1$ 

\medskip

In Case 1 we argue as follows: Let $z\in V$ such that $z^2\not\in yR_1$ (here we use that the characteristic of $K$ is not $2$). Complete $y, z$
to a $K$-basis of $V$ with elements $z_i$. Since $R_2/yR_1$ is $1$-dimensional generated by $z^2$, we may replace $z_i$ with $z_i-*z$ and assume that $z_i^2\in yR_1$. Now, if for some $i$, $z_iR_1\not\subseteq yR_1$ we end up in
case (2) of Lemma \ref{ml}. Hence we have to discuss the case in which $z_iR_1\subseteq yR_1$. In other words, the
$z_i$ are in the socle of $R/(y)$. Modding out this socle elements we get an algebra $S$ with Hilbert series $1+3x+x^2$ and the residue class of $z$ in $S$ satisfies $z^2\neq 0$. So by \ref{l3} there exists $w\in R_1$ such that
$wz\not\in yR_1$ and $w^2\in yR_1$. This is case (3) of \ref{ml}. 
\medskip 

Case 2: Take $z\in V$ such that $zR_1\not\subseteq yR_1$ and note that this is case (2) of
\ref{ml}. 
\medskip 

Case 3: In $R/(y)$ the space $V/(y)$ belongs to the socle. So since $R$ is quadratic and Artinian, the algebra $R/(V)$ has Hilbert series $1+2x+x^2$. For such an algebra it is easy to see that there exists independent linear forms $t,w$ such
 that $t^2=0$ and $w^2=0$. Lifting back to $R$, we have that a basis of $R_2$ is given by $ty, wy, wt$ and for
every $z\in V$ not multiple of
$y$ we get equations: 

\begin{equation}
\label{EQ}
 \quad y^2=0, \ yz=0, \ z^2=L_1y, \ t^2=L_2y, \ w^2=L_3y, \ zw=L_4y, \ zt=L_5y
 \end{equation}

where the $L_i$'s are linear form in $t$ and $w$, say 

$$L_i=\lambda_{i,1}t+\lambda_{i,2}w.$$ Now we look for linear forms of type $\ell=t+a z+b y$ such that $\ell^2=0$. 
The condition $\ell^2=0$ translates into the polynomial system: 

$$\left\{
\begin{array}{ll}
\lambda_{2,1}+2a\lambda_{5,1}+2b+a^2\lambda_{1,1}=0\\
\lambda_{2,2}+2a\lambda_{5,2}+a^2\lambda_{1,2}=0
\end{array}
\right.
$$

Now, assume that 

\begin{equation}
\label{+} 
 \lambda_{1,2}\neq 0 \mbox{ or } \lambda_{5,2}\neq 0.
 \end{equation} 
 
 Then we can solve the second equation to
obtain the value of $a$ and, substituting in the first, we get the value of $b$. In other words, assuming \ref{+}, 
there exists $\ell=t+a z+b y$ such that $\ell^2=0$ in $R$. We evaluate now the rank of such an $\ell$. We have: 

\begin{equation}
\label{B}
\ \ell y= ty \quad \mbox{ and } \quad \ell w=tw+aL_4y+bwy
\end{equation}

and since these two elements of $R_2$ are linearly independent, we can conclude that $\rank
\ell\geq 2$. If $\rank \ell=3$ then we are done by \ref{ml} (1). Hence we may assume that $\rank
\ell=2$. This implies that $\ell z$ and $\ell t$ are linear combinations of $\ell y$ and $\ell w$. Now $\ell
z=(\lambda_{5,2}+a\lambda_{1,2})wy+*ty$ and $\ell t=(\lambda_{2,2}+a\lambda_{5,2})wy+*ty$. Summing up, if $\rank
\ell=2$ then 

\begin{equation}
\label{++}
\quad \lambda_{5,2}+a\lambda_{1,2}=0 \mbox{ and } \lambda_{2,2}+a\lambda_{5,2}=0.
\end{equation} 

In this case, the space
$V(\ell):=\{u\in R_1 : u\ell=0\}$ contains $z+\gamma y$ with
$\gamma=-\lambda_{5,1}-a\lambda_{1,1}$. Note that $(z+\gamma y)^2=\lambda_{1,1}ty+\lambda_{1,2}wy$ and we claim that 
$(z+\gamma y)^2\not\in ell R_1$. If, otherwise, $(z+\gamma y)^2\in \ell R_1$ then, since the $\rank \ell=2$ and hence
the elements \ref{B} are a basis of $R_1$, we must have $\lambda_{1,2}=0$. By \ref{++} follows that $\lambda_{5,2}=0$
condradicting the assumption \ref{+}. We can conclude that $V(\ell)^2 \not\subseteq \ell R_1$. This is Case (1) with $\ell$ playing the role of
$y$ and we are done. Summing up, if \ref{+} holds then we are done because we find an element $\ell$ of with
$\ell^2=0$ which has either rank $3$ or rank $2$ and $\{u\in R_1 : u\ell=0\}^2\not\subseteq \ell R_1$. 

We can aslo look for elements of the form $\ell=w+a z+b y$ satisfying $\ell^2=0$. The situation is completly symmetric. Therefore
we get the desired conclusion unless: 

\begin{equation}
\label{+++} 
\quad \lambda_{1,2}= 0 \mbox{ and } \lambda_{5,2}= 0 \mbox{ and } \lambda_{1,1}=0 \mbox{ and }
\lambda_{4,1}=0
\end{equation}
that is, the equations \ref{EQ} take the form: 
$$ y^2=0, \ \ yz=0, \ \ z^2=0, \ \ t^2=L_2y, \ \ w^2=L_3y, \ \ zw=\lambda_{4,2}wy, \ \ 
zt=\lambda_{5,1}ty$$
 for every $z\in V$ which is not multiple of $y$. 
 It follows that $z^2=0$ for all the $z\in V$. 
This implies $V^2=0$ (here we use again that $K$ has characteristic not $2$). In particular we see that the element
$y_1=z-\lambda_{4,2}y$ of $V$ has $y_1^2=0$ and $y_1V=0$ and $y_1w=0$. So $\rank y_1=1$ and we are done by
\ref{ml} (4). 
\end{proof}

The following examples show that the cases described in \ref{ml} indeed arise in an essential way. By this we
mean, 

\begin{example} \label{ex1} There are examples where all the elements with $\ell^2=0$ have rank $3$. This is the generic
situation. In $4$ variables, an ideal generated by the squares of $7$ general linear forms has this property.
 Explicitely, in $K[t,z,w,y]$ the algebra defined by $(y^2,z^2, w^2, t^2,(t+z+w+y)^2, (t+2z+4w+8y)^2,
(t+3z+9w+27y)^2)$ has this property. 
\end{example}

\begin{example} \label{ex2} There are examples where all the elements with $\ell^2=0$ have rank $2$ and (2) does apply. For example, in $K[t,z,w,y]$ the ideal 
$(y^2,yz,z^2-wy, t^2, tw, w^2-tz,wz)$ defines an algebra $R$ with only two elements with $\ell^2=0$, namely $y$ and
$t$ and both have of rank $2$. For both $y$ and $t$ one can apply \ref{ml} (2). 
\end{example}

\begin{example} \label{ex3} 
There are examples where all the elements with $\ell^2=0$ have rank $2$ and (2) does not apply 
while \ref{ml}(3) does apply. 
The ideal $(y^2,yz, w^2,wz, t^2, tz, z^2+ty+wy)$ defines an algebra with three elements with $\ell^2=0$, namely $y,w,t$. They all have rank $2$. While $y$ does not fit
into \ref{ml} (2) or (3), both $t$ and $w$ satisfy \ref{ml} (3) but not (2).
\end{example}

\begin{example} \label{ex4} There are examples where the elements with $\ell^2=0$ have all rank $\leq 2$ and for those of rank
$2$ case (2) or (3) do not apply. The ideal $(y^2, zy, z^2, t^2-ty-2wy, w^2-3ty-4wy, tz-ty, wz-2wy )$ defines
an algebra where the elements with $\ell^2=0$ are the element of type $ay+bz$. They all have rank $\leq 2$. Those of
rank exactly $2$ do not fit into \ref{ml} (2) or (3). Among the $\ell^2=0$ there are exactly $2$ elements of rank
$1$, namely $y-z$ and $2y-z$. 
 \end{example}
 
 \section{Nets of conics}
 \label{appl} 
Our main result asserts that quadratic Artinian algebras with $\dim R_2=3$ are Koszul and most of them are G-quadratic. 
What about dropping the assumption Artinian? We will discuss in this section the case of quadratic (non-necessarily Artinian)  algebras with Hilbert series $1+3x+3x^2+\dots$.   In \cite{BF} the authors make a detailed study of the Koszul property of the quadratic quotients of $K[x,y,z]$. 
 The most difficult case is that of a quotient defined by $3$ quadrics, that is, algebra with Hilbert series $1+3x+3x^2+\dots$.  
It turns out that there exist exactly (up to change of coordinates) two quotients of $K[x,y,z]$ defined by $3$ quadrics that  are not Koszul. They are the algebras defined by the ideal number 12) and number 14) in the list below. 
 
To proceed with the discussion let us recall few facts.  Vector spaces of quadrics of dimension $3$  in $3$  variables are classically called nets  conics.   The main ingredient for the proof of \cite[Theorem 1]{BF} is a classification result for nets of conics up to the action of $\GL_3(K)$. This classification can be found in full details in the paper of Wall \cite{W} or in a  old preprint of Emsalem and Iarrobino \cite{EI}. Over the complex numbers, there are $15$ types of nets of conics,  fourteen   of them are just one point and one type is $1$-dimensional. With respect \cite{W},  we have chosen slightly different  normal forms  to minimize the total number of terms involved or (as in case 15) to maximize the symmetry. The list of nets of conics is: 
 $$
 \begin{array}{rlrlrl}
1) & (x^2,xy,y^2) & 2) & (x^2,xy,xz) & 3) & (x^2,y^2,z^2) \\
4) & (xy,xz,yz) & 5) & (x^2,y^2,z^2+xy) & 6) & (xz, yz, z^2+xy) \\
7) & (x^2,y^2,xz) & 8) & (xy,z^2,yz) & 9) & (x^2,y^2,xz+yz) \\
10) & (x^2,xz, y^2 + yz) & 11) & (x^2, xy+z^2, xz) & 12) & (x^2, xy+z^2, yz) \\
13) & (x^2, yz, y^2+z^2+xy) & 14) & (xz, y^2+yx, z^2+xy) \\
\end{array} 
$$
and 
$$ 15)\ \  ( x^2+2jyz, y^2+2jxz, z^2+2jxy) \mbox{ with } j\in \CC \mbox{ and } j^3\neq 0,1,-1/8 $$

One can then study the G-quadratic property of the quadratic quotients of $K[x,y,z]$, say over $\CC$,  using the  classification.  We skip the uninteresting details. Below we summarize the final result. 

There are $5$ possible Hilbert series: 
$$
\begin{array}{rlrlrl}
a) & 1+3x+3x^2+x^3 \mbox{ (c.i.) } & b) & (1+2x)/(1-x) \\ 
c) & (1 + x - 2x^2 + x^3)/(1-x)^2  & d) & (1 + 2x - x^3)/(1-x) \\
 e) & (1 + 2x - 2x^3)/(1-x) 
\end{array} 
$$

Every net of conics $V$ has a dual net of conics $V^*$ (the orthogonal space with respect to partial differentiation). In this duality, a point $(a,b,c)\in \PP^2$ belongs to the locus defined by $V$ if and only the square of the linear form $ax+by+cz$ belongs to $V^*$. 

Another interesting aspect of the story is the following. If a net $V$ is generated by the partial derivatives of a cubic $f$ we say that $V$ is of gradient type. It turns out that ``almost all" nets of conics are of gradient type and the corresponding cubic is also uniquely determined. For instance, the net 15) corresponds to the   smooth cubic form in  Hesse  form   $f=x^3+y^3+z^3+6jxyz$ and the conditions on $j$ guarantee that $f$ is smooth and not in the orbit of the Fermat cubic $x^3+y^3+z^3$. 

In the following table we show for every type, its Hilbert series (column-mark H-series), the number of linear forms whose squares are in the net (column-mark q), the number of points of the variety defined by the net (column-mark p), whether it is Koszul or not (column-mark Kos), whether it is G-quadratic or not (column-mark G-quad), the name given by Wall to that type of net (column-mark Wall), and whether it is of gradient type (column-mark $\nabla$). 

\bigskip 

 \centerline{
 \begin{tabular} {|r||c|c|c|c|c|c|c|c|}
 \hline 
 & {H-series} & q & p & {Kos} & {G-quad} & Wall & $\nabla$ 
\\ \hline \hline
1) & b) & $\infty$ & 1 & { yes } & { yes } & I & no\\ 
 \hline 
 2) & c) & 1 & $\infty$ & { yes } & { yes } & I* & no \\ 
 \hline 
3) & a) & 3 & 0 & { yes } & { yes } & E & yes \\
 \hline 
4) & b) & 0 & 3 & { yes }& { yes } & E* & yes\\
 \hline 
5) & a) & 2& 0 & { yes }& { yes } & D & no \\
 \hline 
6) & d) & 0& 2& { yes } & { yes } & D* & yes \\
 \hline 
7) & d) & 2& 1 & { yes } & { yes } & G & yes\\
 \hline 
8) & b) & 1& 2 & { yes } & { yes } & G* & no\\
 \hline 
9) & d) & 2& 1 & { yes } & { yes } & F & no\\
 \hline 
10) & d) & 1& 2 & { yes } & { yes } & F* & no\\
 \hline 
11) & b) & 1& 1& { yes } & { yes } & H & yes\\
 \hline 
12) & e) & 1&1 & { no } & { no } & C & no \\
 \hline 
13) & a) & 1& 0 & { yes } & { yes } & B & no \\
 \hline 
14) & e) & 0 &1 & { no } & { no } & B* & yes \\
 \hline 
 15) & a) & 0 & 0 & { yes } & { no } & A & yes \\
 \hline 
\end{tabular}
}

\bigskip

The star $*$ in Wall's notation refers to the duality. 
 The types that are G-quadratic are so in the given coordinates with the exception of 6). Applying the change of coordinates, $x\to x, y\to x - z, z\to x - y$,  the net  6) becomes  generated by $xz - yz, x^2 - xy, y^2 - yz$ which is a G-basis. That 12) and 14) are not G-quadratic follows from the fact, proved in \cite{BF}, that they are not Koszul. But it follows also from the simple observation that there is no quadratic monomial ideal with  Hilbert series e). 
 
 \section{Final remarks} 
We list some questions which arise from the results presented. 
Let $R$ be a quadratic standard graded $K$-algebra which is not a trivial fiber extension.
\begin{itemize} 
\item[(1)] Assume $\dim R_1>\dim R_2=3$.  Is $R$ G-quadratic? 
\item[(2)] Assume $\dim R_1>\dim R_2$ and $R_3=0$ (or just $R$ is Artinian). Is $R$ G-quadratic? 
\end{itemize} 
 
We have seen that the answer  to (1) is positive for Artinian algebras. Also, the answer to  (2) is positive for  ``generic" algebras, see \cite{C}.

\end{document}